\documentclass[runningheads]{llncs}

\usepackage[utf8]{inputenc}
\usepackage{amsmath}
\usepackage{amssymb}
\usepackage{stmaryrd}
\usepackage{multirow}
\usepackage{graphicx}
\usepackage{dsfont}
\usepackage{url}
\usepackage{subcaption}
\usepackage{soul}
\setuldepth{Paris}

\newcommand{\R}{\mathbb{R}}

\newcommand{\mc}{\mathcal}

\newcommand{\lto}{\longrightarrow}
\newcommand{\lmto}{\longmapsto}

\newcommand{\cns}{\Longleftrightarrow}
\newcommand{\ls}{\leqslant}

\newcommand{\tr}{\mathrm{tr}}

\newcommand{\Diag}{\mathrm{Diag}}

\newcommand{\Sym}{\mathrm{Sym}}
\newcommand{\sym}{\mathrm{S}}
\newcommand{\symm}{\mathrm{sym}}

\newcommand{\GL}{\mathrm{GL}}

\newcommand{\Cor}{\mathrm{Cor}}

\newcommand{\Exp}{\mathrm{Exp}}
\newcommand{\Log}{\mathrm{Log}}

\newcommand{\ver}{\mathrm{ver}}
\newcommand{\hor}{\mathrm{hor}}

\newcommand{\Hol}{\mathrm{Hol}}

\begin{document}

\title{Geodesics of the Quotient-Affine Metrics on Full-Rank Correlation Matrices}
\titlerunning{Geodesics of the Quotient-Affine Metrics}

\author{Yann Thanwerdas \and Xavier Pennec}
\authorrunning{Y. Thanwerdas, X. Pennec}

\institute{Universit\'e C\^ote d'Azur, Inria, Epione, France\\ \mailsa, \mailsb}
\urldef{\mailsa}\path|yann.thanwerdas@inria.fr|
\urldef{\mailsb}\path|xavier.pennec@inria.fr|

\maketitle

\vspace*{-3mm}
\begin{abstract}
Correlation matrices are used in many domains of neurosciences such as fMRI, EEG, MEG. However, statistical analyses often rely on embeddings into a Euclidean space or into Symmetric Positive Definite matrices which do not provide intrinsic tools. The quotient-affine metric was recently introduced as the quotient of the affine-invariant metric on SPD matrices by the action of diagonal matrices. In this work, we provide most of the fundamental Riemannian operations of the quotient-affine metric: the expression of the metric itself, the geodesics with initial tangent vector, the Levi-Civita connection and the curvature.

\keywords{Riemannian geometry \and Quotient manifold \and Quotient-affine metric \and SPD matrices \and Correlation matrices.}
\end{abstract}

\vspace*{-3mm}
\section{Introduction}
\enlargethispage{8mm}
Correlation matrices are used in many domains with time series data such as functional brain connectivity in functional MRI, electroencephalography (EEG) or magnetoencephalography (MEG) signals. However, there are very few tools to compute intrinsically with correlation matrices from a geometric point of view. Full-rank correlation matrices form a strict sub-manifold of the cone of Symmetric Positive Definite (SPD) matrices sometimes called the (open) elliptope \cite{Tropp18}. Since there exist efficient tools on SPD matrices (affine-invariant/Fisher-Rao metric, log-Euclidean metric...), correlation matrices are often treated as SPD matrices \cite{Varoquaux10}. Nevertheless, these extrinsic tools do not respect the geometry of correlation matrices. Moreover, most of these tools on SPD matrices are invariant under orthogonal transformations which is not compatible with correlation matrices. The elliptope is not even stable by the orthogonal action. Hence using the tools designed for SPD matrices may not be relevant for correlation matrices.

Other structures were proposed on correlation matrices. For example, \cite{Rebonato00} rely on a surjection from a product of $n$ spheres of dimension $n-1$ onto the space of correlation matrices in order to sample valid correlation matrices for financial applications: a point in the former space can be represented by an $n\times n$ matrix $A$ with normed rows and therefore encode a correlation matrix $AA^\top$. Since low-rank matrices have null measure, one gets a full-rank correlation matrix almost surely. 
More recently, the open elliptope was endowed with Hilbert's projective geometry \cite{Nielsen19} thanks to its convexity, or with the quotient of the affine-invariant metric by the action of diagonal matrices  \cite{David19,David19-thesis}.

In this work, we investigate the geometry of this last quotient-affine metric and we contribute additional tools to compute on this manifold. 
Based on the formalization of the quotient manifold with vertical and horizontal distributions, we compute in closed form some fundamental Riemannian operations of the quotient-affine metrics, notably the particularly important exponential map.
In Section 2, we recall how quotient-affine metrics are introduced, as well as the basics on quotient manifolds. In Section 3, we provide the following fundamental quotient and Riemannian operations of quotient-affine metrics: the vertical and horizontal projections, the metric, the exponential map, the Levi-Civita connection and the sectional curvature. This opens the way to many practical algorithms on the open elliptope for different applications. Considering SPD matrices as the Cartesian product of positive diagonal matrices and full-rank correlation matrices, it also allows to define new Riemannian metrics which preserves the quotient-affine geometry on correlation matrices. Thus in Section 4, we illustrate the quotient-affine metric in dimension 2 by coupling it with the diagonal power-Euclidean metrics $g^{\mathrm{E}(p)}_D(\Delta,\Delta)=\tr(D^{2(p-1)}\Delta^2)$ for $p\in\{-1,0,1,2\}$ and then by comparing it with the affine-invariant and the log-Euclidean metrics on SPD matrices.

\section{Quotient-affine metrics}
\enlargethispage{5mm}
\subsection{The quotient manifold of full-rank correlation matrices}\label{subsec:quotient_manifold}

The group of positive diagonal matrices $\Diag^+(n)$ acts on the manifold of SPD matrices $\Sym^+(n)$ via the congruence action $(\Sigma,D)\in\Sym^+(n)\times\Diag^+(n)\lmto D\Sigma D\in\Sym^+(n)$.
The manifold of full-rank correlation matrices $\Cor^+(n)$ can  be seen as the quotient manifold $\Sym^+(n)/\Diag^+(n)$ via the invariant submersion $\pi$ which computes the correlation matrix from a covariance matrix, $\pi:\Sigma\in\Sym^+(n)\lmto\Diag(\Sigma)^{-1/2}\,\Sigma\,\Diag(\Sigma)^{-1/2}\in\Cor^+(n)$.
Hence, any Riemannian metric $G$ on $\Sym^+(n)$ which is invariant under $\Diag^+(n)$ induces a quotient metric $g$ on $\Cor^+(n)$. The steps to define it are the following.
\begin{enumerate}
    \item \ul{Vertical distribution}.
    $\mc{V}_\Sigma=\ker d_\Sigma\pi$ for all $\Sigma\in\Sym^+(n)$.
            
    \item \ul{Horizontal distribution}. 
    $\mc{H}_\Sigma:=\mc{V}_\Sigma^\perp$ for all $\Sigma\in\Sym^+(n)$, where the orthogonality $\perp$ refers to the inner product $G_\Sigma$.
            
    \item \ul{Horizontal lift}. The linear map $d_\Sigma\pi$ restricted to the horizontal space $\mc{H}_\Sigma$ is a linear isomorphism onto the tangent space of full-rank correlation matrices $(d_\Sigma\pi)_{|\mc{H}_\Sigma}:\mc{H}_\Sigma\overset\sim\lto T_{\pi(\Sigma)}\Cor^+(n)$. The horizontal lift $\#$ is its inverse:
            \begin{equation}
                \#:X\in T_{\pi(\Sigma)}\Cor^+(n)\overset\sim\lto X^\#\in\mc{H}_\Sigma.
            \end{equation}
    
    \item \ul{Quotient metric}. It is defined by pullback through the horizontal lift:
            \begin{equation}
                \forall C\in\Cor^+(n),\forall X\in T_C\Cor^+(n),g_C(X,X)=G_\Sigma(X^\#,X^\#),
            \end{equation}
            where $\Sigma\in\pi^{-1}(C)$ and the definition does not depend on the chosen $\Sigma$.
\end{enumerate}

So the only missing ingredient is a Riemannian metric on SPD matrices which is invariant under the congruence action of positive diagonal matrices. In \cite{David19,David19-thesis}, the authors chose to use the famous affine-invariant/Fisher-Rao metric.

\subsection{The affine-invariant metrics and the quotient-affine metrics}

The affine-invariant metric is the Riemannian metric defined on SPD matrices by $G_\Sigma(V,V)=\tr(\Sigma^{-1}V\Sigma^{-1}V)$ for all $\Sigma\in\Sym^+(n)$ and $V\in T_\Sigma\Sym^+(n)$ \cite{Siegel43,Skovgaard84,Amari00}. It is invariant under the congruence action of the whole real general linear group $\GL(n)$ which contains $\Diag^+(n)$ as a subgroup. It provides a Riemannian symmetric structure to the manifold of SPD matrices, hence 
it is geodesically complete and the geodesics are given by the one-parameter subgroups and the group action. We recall the exponential map, the Levi-Civita connection and the sectional curvature below for all $\Sigma\in\Sym^+(n)$ and $V,W\in T_\Sigma\Sym^+(n)$:
\small
\begin{align}
    \Exp^G_\Sigma(V)&=\Sigma^{1/2}\exp(\Sigma^{-1/2}V\Sigma^{-1/2})\Sigma^{1/2},\\
    (\nabla^G_VW)_{|\Sigma}&=\partial_VW-\frac{1}{2}(V\Sigma^{-1}W+W\Sigma^{-1}V),\\
    \kappa^G_\Sigma(V,W)&=\frac{\tr(\Sigma^{-1}V\Sigma^{-1}W\Sigma^{-1}(V\Sigma^{-1}W-W\Sigma^{-1}V))}{2\,G(V,V)G(W,W)},
\end{align}
\normalsize
where $\partial$ is the Euclidean connection of $\Sym(n)$ induced on $\Sym^+(n)$.

The metrics that are invariant under the congruence action of the general linear group $\GL(n)$ actually form a two-parameter family of metrics indexed by $\alpha>0$ and $\beta>-\alpha/n$ \cite{Pennec08}: $G^{\alpha,\beta}_\Sigma(V,V)=\alpha\,\tr(\Sigma^{-1}V\Sigma^{-1}V)+\beta\,\tr(\Sigma^{-1}V)^2$.
We call them all affine-invariant metrics. In particular, these metrics are invariant under the congruence action of diagonal matrices so they are good candidates to define Riemannian metrics on full-rank correlation matrices by quotient.
In \cite{David19,David19-thesis}, the authors rely on the ``classical" affine-invariant metric ($\alpha=1$, $\beta=0$). We generalize below their definition.

\begin{definition}[Quotient-affine metrics on full-rank correlation matrices]
The quotient-affine metric of parameters $\alpha>0$ and $\beta>-\alpha/n$ is the quotient metric on $\Cor^+(n)$ induced by the affine-invariant metric $G^{\alpha,\beta}$ via the submersion $\pi:\Sigma\in\Sym^+(n)\lmto\Diag(\Sigma)^{-1/2}\,\Sigma\,\Diag(\Sigma)^{-1/2}\in\Cor^+(n)$.
\end{definition}

\enlargethispage{8mm}

\vspace*{-5mm}
\section{Fundamental Riemannian operations}

In this section, we detail the quotient geometry of quotient-affine metrics. 
We give the vertical and horizontal distributions and projections in Section \ref{subsec:ver_hor}. We contribute the formulae of the metric itself in Section \ref{subsec:metric}, the exponential map in Section \ref{subsec:geodesics}, and finally the Levi-Civita connection and the sectional curvature in Section \ref{subsec:levi-civita_curvature}. To the best of our knowledge, all these formulae  were not given in the Paul David's paper \cite{David19} and thesis \cite{David19-thesis} and they are new.

\subsection{Vertical and horizontal distributions and projections}\label{subsec:ver_hor}

\begin{itemize}
\item Let $\bullet$ be the Hadamard product on matrices defined by $[A\bullet B]_{ij}=A_{ij}B_{ij}$. 
\item Let $A:\Sigma\in\Sym^+(n)\lmto A(\Sigma)=\Sigma\bullet\Sigma^{-1}\in\Sym^+(n)$. This smooth map is invariant under the action of positive diagonal matrices. Note that the Schur product theorem ensures that $\Sigma\bullet\Sigma^{-1}\in\Sym^+(n)$.
\item Let $\psi:\mu\in\R^n\lmto (\mu\mathds{1}^\top+\mathds{1}\mu^\top)\in\Sym(n)$. This is an injective linear map. 
\item Let $\mc{S}_\Sigma(V)$ the unique solution of the Sylvester equation $\Sigma\mc{S}_\Sigma(V)+\mc{S}_\Sigma(V)\Sigma=V$ for $\Sigma\in\Sym^+(n)$ and $V\in\Sym(n)$. 
\item Let $\Hol^\sym(n)$ be the vector space of symmetric matrices with vanishing diagonal (symmetric hollow matrices). Each tangent space of the manifold of full-rank correlation matrices can be seen as a copy of this vector space.
\end{itemize}

\begin{theorem}[Vertical and horizontal distributions and projections]
\label{thm:VerticalProj}
The vertical distribution is given by $\mc{V}_\Sigma=\Sigma\bullet\psi(\R^n)$ and the horizontal distribution is given by $\mc{H}_\Sigma=\mc{S}_{\Sigma^{-1}}(\Hol^\sym(n))$. The vertical projection is:
\begin{equation}
    \ver:V\in T_\Sigma\Sym^+(n)\lmto\Sigma\bullet\psi((I_n+A(\Sigma))^{-1}\Diag(\Sigma^{-1}V)\mathds{1})\in\mc{V}_\Sigma .
\end{equation}
Then, the horizontal projection is simply $\hor(V)=V-\ver(V)$.
\end{theorem}

\subsection{Horizontal lift and metric}\label{subsec:metric}

\begin{theorem}[Horizontal lift]
\label{thm:HorizontalLift}
Let $\Sigma\in\Sym^+(n)$ and $C=\pi(\Sigma)\in\Cor^+(n)$. The horizontal lift at $\Sigma$ of $X\in T_C\Cor^+(n)$ is $X^\#=\hor(\Delta_\Sigma X\Delta_\Sigma)$ with $\Delta_\Sigma=\Diag(\Sigma)^{1/2}$. In particular, the horizontal lift at $C\in\Sym^+(n)$ is $X^\#=\hor(X)$.
\end{theorem}

\begin{theorem}[Expression of quotient-affine metrics]
\label{thm:QAMetric}
The expression of quotient-affine metrics is for all $C\in\Cor^+(n)$ and $X\in T_C\Cor^+(n)$:
\begin{equation}
    g^{\alpha,\beta}_C(X,X)=G^{\alpha,\beta}_C(X,X)-2\mu^\top[\alpha(I_n+A(C))+2\beta\mathds{1}\mathds{1}^\top]\mu.
\end{equation}
where $\mu=(I_n+A(C))^{-1}\Diag(C^{-1}X)\mathds{1}$.
\end{theorem}

\subsection{Geodesics}\label{subsec:geodesics}

The geodesics of a quotient metric are the projections of the horizontal geodesics of the original metric. This allows us to obtain the exponential map of the quotient-affine metrics. 

\begin{theorem}[Geodesics of quotient-affine metrics]
The geodesic from $C\in\Cor^+(n)$ with initial tangent vector $X\in T_C\Cor^+(n)$ is:
\begin{equation}
    \forall t\in\R,\gamma_{(C,X)}(t)=\Exp_C(X)=\pi(C^{1/2}\exp(C^{-1/2}\hor(X)C^{-1/2})C^{1/2}).
\end{equation}
In particular, the quotient-affine metric is geodesically complete.
\end{theorem}

The Riemannian logarithm between $C_1$ and $C_2\in\Cor^+(n)$ is much more complicated to compute. The general method on Riemannian quotient manifolds  consists in finding $\Sigma\in\Sym^+(n)$ in the fiber above $C_2$ such that $\Log^G_{C_1}(\Sigma)$ is horizontal. Then we have $\Log_{C_1}(C_2)=d_{C_1}\pi(\Log^G_{C_1}(\Sigma))$. This means finding $\Sigma$ that minimizes the affine-invariant distance in the fiber: 
\[ D=\mathrm{arg\,min}_{D\in\Diag^+(n)}d(C_1,DC_2D), \] from which we get  $\Sigma=DC_2D$. This is the method used in the original paper \cite{David19,David19-thesis}.
Note that the uniqueness of the minimizer has not been proved yet.

\subsection{Levi-Civita connection and sectional curvature}\label{subsec:levi-civita_curvature}

In this section, we give the Levi-Civita connection and the curvature without detailing the proofs. The computations are based on the fundamental equations of submersions \cite{O'Neill66}. We denote $\symm(M)=\frac{1}{2}(M+M^\top)$ the symmetric part of a matrix.

\begin{theorem}[Levi-Civita connection and sectional curvature of quotient-affine metrics]
The Levi-Civita connection of quotient-affine metrics is:
\small
\begin{align}
    \nonumber& \textstyle (\nabla_XY)_{|C}=(\partial_XY)_{|C}+\symm[\Diag(X^\#)Y^\#+\Diag(Y^\#)X^\#+\Diag(X^\# C^{-1}Y^\#)C\\
    &\textstyle \quad\quad\quad-X^\# CY^\#-\frac{1}{2}\Diag(X^\#)C\Diag(Y^\#)-\frac{3}{2}\Diag(X^\#)\Diag(Y^\#)C].
\end{align}
\normalsize
The curvature of quotient-affine metrics is:
\begin{equation}\textstyle
    \kappa_C(X,Y)=\kappa^G_C(X^\#,Y^\#)+\frac{3}{4}\frac{G_C(\ver[X^\#,Y^\#],\ver[X^\#,Y^\#])}{g_C(X,X)g_C(Y,Y)-g_C(X,Y)^2},
\end{equation}
where $\kappa^G_C(V,W)=\frac{1}{2}\frac{\tr(C^{-1}VC^{-1}WC^{-1}(VC^{-1}W+WC^{-1}V))}{G_C(V,V)G_C(W,W)-g_C(V,W)^2}\ls 0$ is the sectional curvature of the affine-invariant metrics,
$[V,W]=\partial_VW-\partial_WV$ is the Lie bracket on $\Sym^+(n)$ and:
\small
\[
    (\partial_{X^\#}Y^\#)_{|C} =\frac{1}{2}(\Diag(X^\#)Y+Y\Diag(X^\#)) 
    -2X^\#\bullet\psi((I_n+A(C))^{-1}\Diag(C^{-1}Y)\mathds{1})+v_C(X,Y),
\]
\normalsize
where $v_C(X,Y)= \textstyle 2C\bullet\psi((I_n+A(C))^{-1}D(X,Y)\mathds{1})\in\mc{V}_C$ and
\small 
\begin{align*}
    D(X,Y)&=\textstyle \Diag(C^{-1}XC^{-1}Y)\\
    &\textstyle \quad+(X^\#\bullet C^{-1}-C\bullet(C^{-1}X^\# C^{-1}))(I_n+A(C))^{-1}\Diag(C^{-1}Y).
\end{align*}
\end{theorem}
\normalsize

Note that the first term of the sectional curvature is negative and the second one is positive so we don't know in general the sign of the curvature of quotient-affine metrics. Although these formulae are quite heavy, it is possible to implement them. This will be the goal of a future work.

The quotient-affine metrics not only provide a Riemannian framework on correlation matrices but also provide correlation-compatible statistical tools on SPD matrices if we consider that the space of SPD matrices is the Cartesian product of positive diagonal matrices and full-rank correlation matrices. We give a taste of such construction in the next section.

\section{Illustration in dimension 2}

In dimension 2, any correlation matrix $C\in\Cor^+(2)$ writes $C=C(\rho):=\begin{pmatrix}1 & \rho\\\rho & 1\end{pmatrix}$ with $\rho\in(-1,1)$. In the following theorem, we give explicitly the logarithm, the distance and the interpolating geodesics of the quotient-affine metric.

\begin{theorem}[Quotient-affine metrics in dimension 2]
Let $C_1=C(\rho_1),C_2=C(\rho_2)\in\Cor^+(n)$ with $\rho_1,\rho_2\in(-1,1)$. We denote $f:\rho\in(-1,1)\mapsto\frac{1+\rho}{1-\rho}\in(0,\infty)$ which is a smooth increasing map. We denote $\lambda=\lambda(\rho_1,\rho_2)=\frac{1}{2}\log\frac{f(\rho_2)}{f(\rho_1)}$ which has the same sign as $\rho_2-\rho_1$. Then:
\begin{enumerate}
    \item $\Log_{C_1}(C_2)=\lambda\begin{pmatrix}0 & 1-\rho_1^2\\1-\rho_1^2 & 0\end{pmatrix}$ and $d(C_1,C_2)=\sqrt{2}|\lambda|$,
    \item $\gamma^{\mathrm{QA}}_{(C_1,C_2)}(t)=C(\rho^{\mathrm{QA}}(t))$ where $\rho^{\mathrm{QA}}(t)=\frac{\rho_1\cosh(\lambda t)+\sinh(\lambda t)}{\rho_1\sinh(\lambda t)+\cosh(\lambda t)}\in(-1,1)$ is monotonic (increasing if and only if $\rho_2-\rho_1>0$).
\end{enumerate}
\end{theorem}

\begin{figure}[h!]
    \centering
    \begin{minipage}{.54\textwidth}
    \subcaptionbox{Interpolation and extrapolation of the correlation coefficient. All the interpolations relying on product metrics on $\Sym^+(2)=\Diag^+(2)\times\Cor^+(2)$ with the quotient-affine metric on $\Cor^+(2)$ lead to the same correlation coefficient, labelled as ``Quotient-affine".}{
    \includegraphics[width=\textwidth]{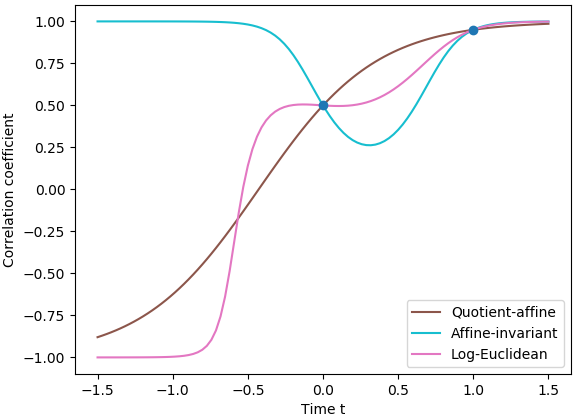}
    \label{fig:curves:cor}}
    \end{minipage}
    \begin{minipage}{.42\textwidth}
    \subcaptionbox{Interpolation of the determinant}{
    \includegraphics[width=\textwidth]{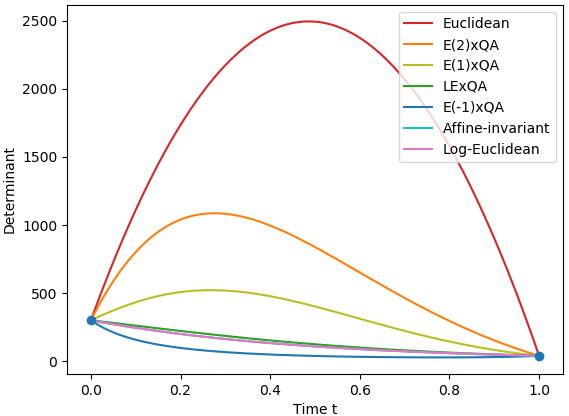}
    \label{fig:curves:det}}
    \subcaptionbox{Interpolation of the trace}{
    \includegraphics[width=\textwidth]{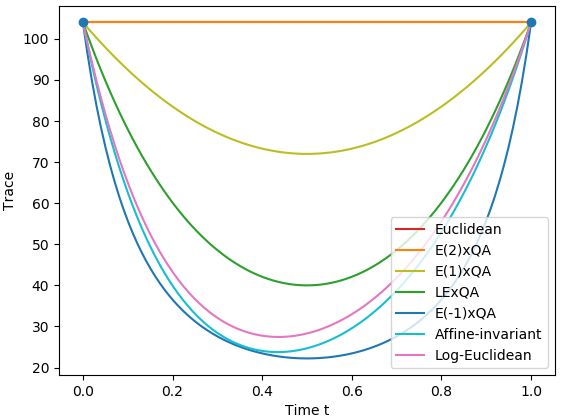}
    \label{fig:curves:tr}}
    \end{minipage}
    \caption{Extrapolation and interpolation between the SPD matrices $\Sigma_1$ and $\Sigma_2$ using various Riemannian metrics. The metrics $\mathrm{E}(p)\mathrm{xQA}$ refer to the $p$-power-Euclidean metric on the diagonal part and the quotient-affine metric on the correlation part. When $p$ tends to $0$, the $p$-power-Euclidean metric $\mathrm{E}(p)$ tends to the log-Euclidean metric $\mathrm{LE}=\mathrm{E}(0)$.}
    \label{fig:curves}
\end{figure}

Let $\Sigma_1,\Sigma_2\in\Sym^+(2)$ and $C_1,C_2$ their respective correlation matrices. We denote $\gamma^{\mathrm{AI}},\gamma^{\mathrm{LE}}$ the geodesics between $\Sigma_1$ and $\Sigma_2$ for the affine-invariant and the log-Euclidean metrics respectively. We define $\rho^{\mathrm{AI}},\rho^{\mathrm{LE}}$ such that the correlation matrices of $\gamma^{\mathrm{AI}}(t),\gamma^{\mathrm{LE}}(t)$ are $C(\rho^{\mathrm{AI}}(t)),C(\rho^{\mathrm{LE}}(t))$. Figure \ref{fig:curves}.(a) shows $\rho^{\mathrm{AI}},\rho^{\mathrm{LE}},\rho^{\mathrm{QA}}$ with $\Sigma_1=\begin{pmatrix}4 & 1\\ 1 & 100\end{pmatrix}$ and $\Sigma_2=\begin{pmatrix}100 & 19\\19 & 4\end{pmatrix}$.
When varying numerically $\Sigma_1$ and $\Sigma_2$, it seems that $\rho^{\mathrm{LE}}$ and $\rho^{\mathrm{AI}}$ always have three inflection points. In contrast, $\rho^{\mathrm{QA}}$ always has one inflection point since $(\rho^{\mathrm{QA}})''=-2\lambda^2\rho^{\mathrm{QA}}(1-(\rho^{\mathrm{QA}})^2)$. Analogously, we compare the interpolations of the determinant (Fig. \ref{fig:curves}.(b)) and the trace (Fig. \ref{fig:curves}.(c)) using several Riemannian metrics: Euclidean (trace-monotonic); log-Euclidean and affine-invariant (determinant-monotonic); power-Euclidean x quotient-affine (correlation-monotonic).

\begin{figure}[t]
    \centering
    \includegraphics[width=\textwidth]{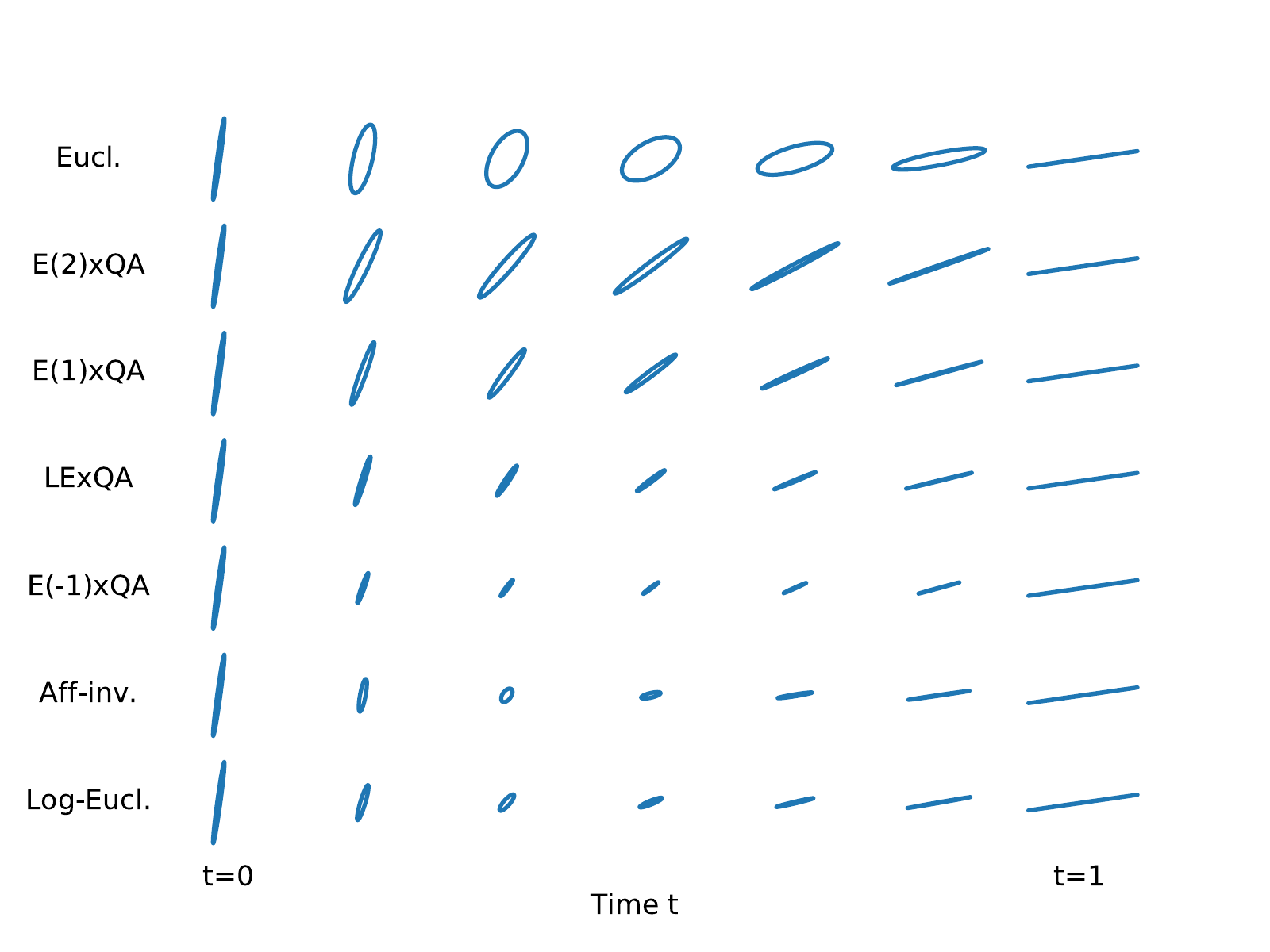}
    \caption{Interpolations between SPD matrices $\Sigma_1$ and $\Sigma_2$.}
    \label{fig:ellipses}
\end{figure}

On Figure \ref{fig:ellipses}, we compare the geodesics of these metrics. The determinant is the area of the ellipsis and the trace is the sum of the lengths of the axes. Thus the product metrics of the form power-Euclidean on the diagonal part and quotient-affine on the correlation part can be seen as performing a correlation-monotonic trade-off between the trace-monotonicity and the determinant-monotonicity.

\section{Conclusion}

We investigated in this paper the very nice idea of quotienting the affine-invariant metrics on SPD matrices by the action of the positive diagonal  group to obtain the principled quotient-affine metrics on full-rank correlation matrices. The quotient-affine metric with $\alpha=1$ and $\beta=0$ was first proposed in Paul David's  thesis \cite{David19-thesis} and in the subsequent journal paper \cite{David19}. We contribute here exact formulae for  the main Riemannian operations, including the exponential map, the connection and the sectional curvature. The exponential map is particularly interesting to rigorously project tangent computations to the space of full-rank correlation matrices. This opens the way to the implementation of a number of generic algorithms on manifolds. 
However, we could not find a closed form expression for the logarithm which remains to be computed through an optimization procedure. Thus, computing distances with these metrics remains computationally expensive.  
In order to obtain more efficient methods, this leads us to look for other principled Riemannian metrics on correlation matrices for which the logarithm could be expressed in closed form.

\textit{Acknowledgements.} This project has received funding from the European Research Council (ERC) under the European Union’s Horizon 2020 research and innovation program (grant G-Statistics agreement No 786854). This work has been supported by the French government, through the UCAJEDI Investments in the Future project managed by the National Research Agency (ANR) with the reference number ANR-15-IDEX-01 and through the 3IA Côte d’Azur Investments in the Future project managed by the National Research Agency (ANR) with the reference number ANR-19-P3IA-0002. The authors would like to thank Frank Nielsen for insightful discussions on correlation matrices.

\bibliographystyle{unsrt}
\bibliography{biblio}

\newpage
\appendix

\small

\section{Proof of theorem \ref{thm:VerticalProj}}

We denote $\Delta_\Sigma=\Diag(\Sigma)^{1/2}$. Using Section \ref{subsec:quotient_manifold}, we have: 
\begin{align*}
    d_\Sigma\pi(V)  &= \textstyle \Delta_\Sigma^{-1}\left[V-\frac{1}{2}(\Delta_\Sigma^{-2}\Diag(V)\Sigma+\Sigma\Diag(V)\Delta_\Sigma^{-2})\right]\Delta_\Sigma^{-1}.\\
    V\in\mc{V}_\Sigma & \textstyle \cns \forall i,j\in\{1,...,n\},V_{ij} =  \Sigma_{ij}\times\frac{1}{2}\left(\frac{V_{ii}}{\Sigma_{ii}}+\frac{V_{jj}}{\Sigma_{jj}}\right)\cns V\in\Sigma\bullet\psi(\R^n).\\
    W\in\mc{H}_\Sigma & \textstyle \cns\forall V\in\mc{V}_\Sigma,\tr(\Sigma^{-1}V\Sigma^{-1}W)=0,\\
    &\textstyle \cns\forall\mu\in\R^n,\sum_{ijkl}{[\Sigma^{-1}]_{ij}\Sigma_{jk}(\mu_j+\mu_k)[\Sigma^{-1}]_{kl}W_{li}}=0,\\
    &\textstyle \cns\forall D\in\Diag(n),\tr(\Sigma^{-1}DW+D\Sigma^{-1}W)=0,\\
    &\textstyle \cns(\Sigma^{-1}W+W\Sigma^{-1})\in\Hol^\sym(n)\cns W\in\mc{S}_{\Sigma^{-1}}(\Hol^\sym(n)).
\end{align*}
Thus, we have computed the horizontal space for the affine-invariant metric $\alpha=1$ and $\beta=0$. It is still valid for all $\alpha>0$ and $\beta>-\alpha/n$ since the latter is included in the former (because $\tr(\Sigma^{-1}W)=0$) and they have the same dimension so they are equal. Now we compute the vertical projection. Let $V\in T_\Sigma\Sym^+(n)$ and let $\mu\in\R^n$ and $W\in\mc{H}_\Sigma$ such that $V=\Sigma\bullet\psi(\mu)+W$. We are looking for $\mu$. Since $\Sigma^{-1}V=\Sigma^{-1}(\Sigma\bullet\psi(\mu))+\Sigma^{-1}W$, we have:
\begin{align*}
    [\Sigma^{-1}V]_{ii}&= \textstyle \sum_j{[\Sigma^{-1}]_{ij}\Sigma_{ij}(\mu_i+\mu_j)}=\mu_i+[A(\Sigma)\mu]_i.\\
    \Diag(\Sigma^{-1}V)\mathds{1}&=(I_n+A(\Sigma))\mu,\mathrm{~hence~}\mu=(I_n+A(\Sigma))^{-1}\Diag(\Sigma^{-1}V)\mathds{1}.
\end{align*}

\section{Proof of theorem \ref{thm:HorizontalLift}}

We compute $d_\Sigma\pi(V)$ with $V=\Delta_\Sigma X\Delta_\Sigma$. Note that $\Diag(V)=0$.
\begin{equation*} \textstyle
    d_\Sigma\pi(V)=\Delta_\Sigma^{-1}\left[V-\frac{1}{2}(\Delta_\Sigma^{-2}\Diag_\Sigma(V)\Sigma+\Sigma\Diag(V)\Delta_\Sigma^{-2})\right]\Delta_\Sigma^{-1}=X.
\end{equation*}
Hence, $\hor(V)\in\mc{H}_\Sigma$ satisfies $d_\Sigma\pi(\hor(V))=X$ so $X^\#=\hor(V)$.

\section{Proof of theorem \ref{thm:QAMetric}}

We use the definition of the quotient metric.
\begin{equation*}
    g^{\alpha,\beta}_C(X,X)=G^{\alpha,\beta}_C(\hor(X),\hor(X))=G^{\alpha,\beta}_C(X,X)-G^{\alpha,\beta}_C(\ver(X),\ver(X)).
\end{equation*}
We denote $\ver(X)=C\bullet\psi(\mu)$ with $\mu=(I_n+A(C))^{-1}\Diag(C^{-1}X)\mathds{1}$. Then:
\begin{align*}
    G_C(\ver(X),\ver(X)) &=\textstyle \,\alpha\,\underset{\sum_{ijkl}{[C^{-1}]_{ij}C_{jk}(\mu_j+\mu_k)[C^{-1}]_{kl}C_{li}(\mu_l+\mu_i)}}{\underbrace{\tr(C^{-1}\ver(X)C^{-1}\ver(X))}}+\beta\,\underset{(\sum_{ij}{[C^{-1}]_{ij}C_{ij}(\mu_i+\mu_j)})^2}{\underbrace{\tr(C^{-1}\ver(X))^2}}\\
    &=\textstyle \,2\alpha\sum_{ij}{(\delta_{ij}+[C^{-1}]_{ij}C_{ij})\mu_i\mu_j}+4\beta\,(\mu^\top\mathds{1})^2\\
    &=\textstyle \,2\mu^\top[\alpha(I_n+A(C))+2\beta\mathds{1}\mathds{1}^\top]\mu.
\end{align*}

\end{document}